\documentclass[12pt]{article}
\usepackage{amssymb, amsmath, amsfonts, url, graphicx}

\def\0{\leqno}
\def\({\left(}
\def\){\right)}
\def\<{\left<}
\def\>{\right>}

\def\4{\subseteq }
\def\dd{\displaystyle}

\def\bit{\begin{itemize}}
\def\eit{\end{itemize}}
\def\barr{\begin{array}}
\def\earr{\end{array}}

\def\X#1#2{\stb{#1}{#2}{\mbox{\Huge$\times$}}}
\def\bld#1#2{{\buildrel{#1}\over{#2}}}
\def\st#1#2{{\mathrel{\mathop{#2}\limits_{#1}}{}\!}}
\def\stb#1#2#3{{\st{{#1}}{\bld{{#2}}{#3}}{}\!}}
\def\xmare#1#2{\stb{#1}{#2}{\mbox{\Huge$\times$}}}
\def\a{{\alpha}}
\def\b{{\beta}}

\def\dd{\displaystyle}

\title{\bf On the sum of element orders\\ of finite abelian groups}
\author{Marius T\u arn\u auceanu, Dan Gregorian Fodor}
\date{October 1, 2014}

\begin{document}

\maketitle

\begin{abstract}
In this note some properties of the sum of element orders of a
finite abelian group are studied.
\end{abstract}



\section{Introduction}

Let $G$ be a finite group. We define the function
$$\psi(G)=\dd\sum_{a\in G}o(a),$$where $o(a)$ denotes the order of $a\in G$.
The starting point for our discussion is given by the papers
\cite{1,2} which investigate the minimum/maximum of $\psi$ on the
groups of the same order.

Recall that the function $\psi$ is multiplicative, that is if
$G_1$ and $G_2$ are two finite groups satisfying $\gcd(|G_1|,
|G_2|)=1$, then
$\psi(G_1\hspace{-1mm}\times\hspace{-1mm}G_2)=\psi(G_1)\psi(G_2)$.
By a standard induction argument, it follows that if $G_{i}$,
$i=1,2,\dots,k$, are finite groups of coprime orders, then
$$\psi(\xmare{i=1}k G_i)=\prod_{i=1}^{k} \psi(G_i)\,.$$This shows
that the study of $\psi(G)$ for finite nilpotent groups $G$ can be
reduced to $p$-groups.

In the current note we will focus on the restriction of $\psi$ to
the class of finite abelian groups $G$. In this case we are able
to give an explicit formula for $\psi(G)$. We prove that abelian
$p$-groups of a fixed order are determined by this quantity and we
conjecture that this happens for \textit{arbitrary} finite abelian
groups. Other interesting properties of the function $\psi$ will
be also discussed.

Most of our notation is standard and will not be repeated here.
Basic concepts and results on group theory can be found in
\cite{3,4}. For subgroup lattice notions we refer the reader to
\cite{5}.

\section{Main results}

As we have seen above, computing the sum of element orders of
finite abelian groups is reduced to $p$-groups. For such a group
$G$ we can determine $\psi(G)$ by using Corollary 4.4 of \cite{6}.

\bigskip\noindent{\bf Theorem 1.} {\it Let $G=\X{i=1}k\,\mathbb{Z}_{p^{\a_i}}$ be a finite abelian
$p$-group, where $1\leq\a_1\leq\a_2\leq...\leq\a_k$. Then
$$\psi(G)=1+\dd\sum_{\a=1}^{\a_k}\left(p^{2\a}f_{(\a_1,\a_2,...,\a_k)}(\a)-p^{2\a-1}f_{(\a_1,\a_2,...,\a_k)}(\a-1)\right),$$where
$$f_{(\a_1,\a_2,...,\a_k)}(\a)=\left\{\barr{lll}
p^{(k-1)\a},&\mbox{ if }&0\le\a\le\a_1\vspace*{1,5mm}\\
p^{(k-2)\a+\a_1},&\mbox{ if }&\a_1\le\a\le\a_2\\
\vdots\\
p^{\a_1+\a_2+...+\a_{k-1}},&\mbox{ if
}&\hspace{-1mm}\a_{k-1}\le\a\,.\earr\right.$$}

\bigskip\noindent{\bf Remarks.}
\begin{itemize}
\item[1.] The function $f_{(\a_1,\a_2,...,\a_k)}$ in Theorem 1 is increasing.
\item[2.] $\psi(\X{i=1}k\,\mathbb{Z}_{p^{\a_i}})$ is a polynomial in $p$ of degree $2\a_k+\a_{k-1}+...+\a_1$.
\item[3.] An alternative way to write $\psi(\X{i=1}k\,\mathbb{Z}_{p^{\a_i}})$
is
$$\psi(\X{i=1}k\,\mathbb{Z}_{p^{\a_i}})=p^{2\a_k+\a_{k-1}+...+\a_1}-\left(p-1\right)\dd\sum_{\a=0}^{\a_k-1}p^{2\a}f_{(\a_1,\a_2,...,\a_k)}(\a).$$
\end{itemize}
\bigskip

Theorem 1 allows us to obtain a precise expression of $\psi(G)$
for some particular finite abelian $p$-groups $G$.

\bigskip\noindent{\bf Corollary 2.} {\it We have:
\begin{itemize}
\item[{\rm a)}] $\psi(\mathbb{Z}_{p^n})=\dd\frac{p^{2n+1}+1}{p+1}$;
\item[{\rm b)}] $\psi(\mathbb{Z}_p^n)=p^{n+1}-p+1$;
\item[{\rm c)}] $\psi(\mathbb{Z}_{p^2}\hspace{-1mm}\times\hspace{-1mm}\mathbb{Z}_p^{n-2})=p^{n+2}-p^{n+1}+p^n-p+1$;
\item[{\rm d)}] $\psi(\mathbb{Z}_{p^{\a_1}}\hspace{-1mm}\times\hspace{-1mm}\mathbb{Z}_{p^{\a_2}})=\dd\frac{p^{2\a_2+\a_1+3}+p^{2\a_2+\a_1+2}+p^{2\a_2+\a_1+1}+p^{3\a_1+2}+p+1}{(p+1)(p^2+p+1)}\,$;
\item[{\rm e)}] $\psi(\mathbb{Z}_{p^{\a_1}}\hspace{-1mm}\times\hspace{-1mm}\mathbb{Z}_{p^{\a_2}}\hspace{-1mm}\times\hspace{-1mm}\mathbb{Z}_{p^{\a_3}})=\dd\frac{p^{2\a_3+\a_2+\a_1+1}+p^{3\a_2+\a_1+2}}{p+1}-\dd\frac{p^{3\a_2+\a_1+3}-p^{4\a_1+3}}{p^2+p+1}-\dd\frac{p^{4\a_1+4}-1}{p^3+p^2+p+1}\,$.
\end{itemize}}
\bigskip

Given a positive integer $n$, it is well-known that there is a
bijection between the set of types of abelian groups of order
$p^n$ and the set $P_n=\{(x_1,x_2,...,x_n)\in\mathbb{N}^n\mid
x_1\geq x_2\geq...\geq x_n, x_1+x_2+...+x_n=n\}$ of partitions of
$n$, namely the map
$$\X{i=1}k\,\mathbb{Z}_{p^{\a_i}} (\mbox{ with }\a_1\leq\a_2\leq...\leq\a_k \mbox{ and }\sum_{i=1}^k \a_i=n)\longmapsto (\a_k,...,\a_1,\hspace{-2mm}\underbrace{0,...,0}_{n-k\,
\rm positions }\hspace{-2mm})\,.$$Moreover, recall that $P_n$ is
totally ordered under the lexicographic order $\preceq$\,, where
$$(x_1,x_2,...,x_n)\hspace{-1mm}\prec\hspace{-1mm}(y_1,y_2,...,y_n)\hspace{-1mm}\Longleftrightarrow\hspace{-1mm} \left\{\barr{lll}
x_1=y_1,...,x_m=y_m\\
\mbox{and}\\
x_{m+1}{<} y_{m+1} \mbox{ for some }
m\hspace{-1mm}\in\hspace{-1mm}\left\{0,1,...,n{-}1\right\}.\earr\right.$$Obviously,
the lexicographic order induces a total order on the set of types
of abelian $p$-groups of order $p^n$.
\bigskip

By computing the values of $\psi$ corresponding to all types of
abelian $p$-groups of order $p^2$, $p^3$ and $p^4$, respectively,
one obtains:
\begin{itemize}
\item[-] $\psi(\mathbb{Z}_p^2)=p^3-p+1<\psi(\mathbb{Z}_{p^2})=p^4-p^3+p^2-p+1$;
\item[-] $\psi(\mathbb{Z}_p^3)=p^4-p+1<\psi(\mathbb{Z}_p\hspace{-1mm}\times\hspace{-1mm}\mathbb{Z}_{p^2})=p^5-p^4+p^3-p+1<\psi(\mathbb{Z}_{p^3})\hspace{-1mm}=p^6-p^5+p^4-p^3+p^2-p+1$;
\item[-] $\psi(\mathbb{Z}_p^4)=p^5-p+1<\psi(\mathbb{Z}_p\hspace{-1mm}\times\hspace{-1mm}\mathbb{Z}_p\hspace{-1mm}\times\hspace{-1mm}\mathbb{Z}_{p^2})=p^6-p^5+p^4-p+1<\psi(\mathbb{Z}_{p^2}\hspace{-1mm}\times\hspace{-1mm}\mathbb{Z}_{p^2})=p^6-p^4+p^3-p+1<\psi(\mathbb{Z}_p\hspace{-1mm}\times\hspace{-1mm}\mathbb{Z}_{p^3})\hspace{-1mm}=p^7-p^6+p^5-p^4+p^3-p+1<\psi(\mathbb{Z}_{p^4})=p^8-p^7+p^6-p^5+p^4-p^3+p^2-p+1$.
\end{itemize}
The above inequalities suggest us that the function $\psi$ is
strictly increasing. This is true, as shows the following theorem.

\bigskip\noindent{\bf Theorem 3.} {\it Let $G_1=\X{i=1}k\,\mathbb{Z}_{p^{\a_i}}$ and $G_2=\X{j=1}r\,\mathbb{Z}_{p^{\b_j}}$ be two finite abelian
$p$-groups of order $p^n$. Then
$$\psi(G_1)<\psi(G_2)\Longleftrightarrow(\a_k,...,\a_1,\hspace{-2mm}\underbrace{0,...,0}_{n-k\,
\rm positions
}\hspace{-2mm})\prec(\b_r,...,\b_1,\hspace{-2mm}\underbrace{0,...,0}_{n-r\,
\rm positions }\hspace{-2mm})\,.\0(*)$$}

\noindent{\bf Proof.} First of all, we remark that it suffices to
prove $(*)$ only for consecutive partitions of $n$ because $P_n$
is fully ordered.

Assume that $(\a_k,...,\a_1,0,...,0)\prec(\b_r,...,\b_1,0,...,0)$.
We have to prove $\psi(G_1)<\psi(G_2)$ (notice that this
inequality holds for the first two elements of $P_n$, by b) and c)
of Corollary 2). Let $s\in\{1,2,...,r-1\}$ such that
$\b_1=\b_2=\cdots=\b_s<\b_{s+1}$. We distinguish the following two
cases.

\medskip\noindent{\bf \hspace{20mm}Case 1.} $\b_1\geq 2$

\noindent Then $(\a_k,...,\a_1,0,...,0)$ is of type
$(\b_r,...,\b_2,\b_1-1,1,0,...,0)$, i.e. $k=r+1$, $\a_1=1$,
$\a_2=\b_1-1$ and $\a_i=\b_{i-1}$ for $i=3,4,...,r+1$. We infer
that
$f_{(\a_1,\a_2,...,\a_k)}(\gamma)=f_{(\b_1,\b_2,...,\b_r)}(\gamma),
\forall \hspace{1mm}\gamma\geq\b_1$. One obtains
$$\psi(G_2)-\psi(G_1)=$$
$$=p^{\b_r+n}-\left(p-1\right)\dd\sum_{\gamma=0}^{\b_r-1}p^{2\gamma}f_{(\b_1,\b_2,...,\b_r)}(\gamma)-p^{\a_k+n}+\left(p-1\right)\dd\sum_{\gamma=0}^{\a_k-1}p^{2\gamma}f_{(\a_1,\a_2,...,\a_k)}(\gamma)=$$
$$\hspace{-42,5mm}=\left(p-1\right)\dd\sum_{\gamma=1}^{\b_1-1}p^{2\gamma}\left(f_{(\a_1,\a_2,...,\a_k)}(\gamma)-f_{(\b_1,\b_2,...,\b_r)}(\gamma)\right)=$$
$$\hspace{-61mm}=\left(p-1\right)\dd\sum_{\gamma=1}^{\b_1-1}p^{2\gamma}\left(p^{(r-1)\gamma+1}-p^{(r-1)\gamma}\right)>0.$$

\medskip\noindent{\bf \hspace{20mm}Case 2.} $\b_1=1$

\noindent Then $(\a_k,...,\a_1,0,...,0)$ is of type
$(\b_r,...,\b_{s+1}-1,\b'_t,\b'_{t-1},...,\b'_1,0,...,0)$, where
$\b_{s+1}-1\geq\b'_t\geq\b'_{t-1}\geq...\geq\b'_1\geq 1$ and
$\b'_t+\b'_{t-1}+...+\b'_1=s+1$. We infer that
$f_{(\a_1,\a_2,...,\a_k)}(\gamma)=f_{(\b_1,\b_2,...,\b_r)}(\gamma),
\forall \hspace{1mm}\gamma\geq\b_{s+1}$. So, we can suppose that
$s=r-1$, i.e.
$$(\a_k,...,\a_1,0,...,0)=(\b_{r}-1,\b'_t,\b'_{t-1},...,\b'_1,0,...,0).$$One
obtains
$$\psi(G_2)-\psi(G_1)=p^{\b_r+n}-\left(p-1\right)\dd\sum_{\gamma=0}^{\b_r-1}p^{2\gamma}f_{(\b_1,\b_2,...,\b_r)}(\gamma)-p^{\b_r+n-1}+S,$$where
$$S=\left(p-1\right)\dd\sum_{\gamma=0}^{\b_r-2}p^{2\gamma}f_{(\a_1,\a_2,...,\a_k)}(\gamma)>
0.$$Since
$$f_{(\b_1,\b_2,...,\b_r)}(\gamma)=\left\{\barr{lll}
p^{(r-1)\gamma},&\mbox{ if }&0\le\gamma\le 1\vspace*{1,5mm}\\
p^{r-1},&\mbox{ if }&1\le\gamma\,,\earr\right.$$it follows that
$$\psi(G_2)-\psi(G_1)> p^{\b_r+n}-p^{\b_r+n-1}-\left(p-1\right)\dd\sum_{\gamma=0}^{\b_r-1}p^{2\gamma}f_{(\b_1,\b_2,...,\b_r)}(\gamma)=$$
$$\hspace{25mm}=p^{\b_r+n}-p^{\b_r+n-1}-\left(p-1\right)\left[1+p^{r-1}\dd\frac{p^{2\b_r}-p^2}{p^2-1}\right]=$$
$$\hspace{26,5mm}=\dd\frac{1}{p+1}\left[p^{\b_r+n-1}\left(p^2-p-1\right)+p^{r+1}-p^2+1\right]>0.$$
\bigskip

Conversely, assume that $\psi(G_1)<\psi(G_2)$, but
$(\a_k,...,\a_1,0,...,0)\succeq(\b_r,...,\newline\b_1,0,...,0)$.
Then the first part of the proof leads to
$\psi(G_2)\leq\psi(G_1)$, a contradiction. Hence $(*)$ holds.
$\scriptstyle\Box$\bigskip

Two immediate consequences of Theorem 3 are the following.

\bigskip\noindent{\bf Corollary 4.} {\it Let $n$ be a positive integer and $G$ be an abelian $p$-group of order
$p^n$. Then the minimum value of $\psi(G)$ is obtained for $G$
elementary abelian, while the maximum value of $\psi(G)$ is obtained
for $G$ cyclic.}

\bigskip\noindent{\bf Corollary 5.} {\it Two finite abelian $p$-groups of the same order are isomorphic if and only if
they have the same sum of element orders.}
\bigskip

Inspired by Corollary 5, we came up with the following conjecture,
which we have verified by computer for all abelian groups of order
less or equal to 100000.

\bigskip\noindent{\bf Conjecture 6.} {\it Two finite abelian groups of the same order are isomorphic if and only if
they have the same sum of element orders.}
\bigskip

In order to decide if two finite abelian groups $G_1$ and $G_2$
are isomorphic by using the above results, the condition $\mid
G_1\hspace{-1mm}\mid\,=\mid G_2\hspace{-1mm}\mid$ is essential, as
shows the following simple example.

\bigskip\noindent{\bf Example.} We have $\mathbb{Z}_2^2\not\cong\mathbb{Z}_3$ even if $\psi(\mathbb{Z}_2^2)=\psi(\mathbb{Z}_3)=7$.
\bigskip

This proves that the function $\psi$ is not injective, too. The
surjectivity of $\psi$ also fails because $Im(\psi)$ contains only
odd positive integers (notice that in fact more can be said,
namely: $\psi(G)$ \textit{is odd for all finite groups} $G$).
Moreover, there exist odd positive integers not contained in
$Im(\psi)$, as 5.
\bigskip

Finally, we observe that $\psi(G)$ is not divisible by
$\mid\hspace{-1mm} G\hspace{-1mm}\mid$ for large classes of finite
groups $G$, as $p$-groups, groups of order $p^nq$ ($p,q$ primes)
without normal Sylow $q$-subgroups and groups of even order. More
precisely, by MAGMA we checked that there are only three types of
groups of order at most 2000 satisfying $\mid
G\hspace{-1mm}\mid\hspace{2mm}\mid \psi(G)$ (the smallest order of
such a group $G$ is 105 and $\psi(G)=1785=105\cdot17$) and these
are not abelian. Consequently, the study of this property for
abelian groups seems to be interesting.

\bigskip\noindent{\bf Theorem 7.} {\it There are finite abelian groups $G$
such that $$\psi(G)\equiv 0\hspace{1mm}({\rm mod} \mid
G\hspace{-1mm}\mid).$$}

\noindent{\bf Proof.} Let
$G=\mathbb{Z}_{13}\hspace{-1mm}\times\hspace{-1mm}\mathbb{Z}_{13}\hspace{-1mm}\times\hspace{-1mm}\mathbb{Z}_{23}$.
We have $\mid G\hspace{-1mm}\mid\hspace{1mm}=3887$ and
$$\psi(G)=\psi(\mathbb{Z}_{13}\hspace{-1mm}\times\hspace{-1mm}\mathbb{Z}_{13})\psi(\mathbb{Z}_{23})=\left(13^3-13+1\right)\dd\frac{23^3+1}{23+1}=$$
$$\hspace{-34mm}=1107795=3887\cdot285,$$completing the proof.
$\scriptstyle\Box$\bigskip

We end our note by indicating a natural generalization of
$\psi(G)$, which is obtained by replacing the orders of elements
with the orders of elements relative to a certain subgroup of $G$.

\bigskip\noindent{\bf Open problem.} Let $G$ be a finite group.
For every subgroup $H$ of $G$, we define the function
$$\psi_H(G)=\dd\sum_{a\in G}o_H(a),$$where $o_H(a)$ denotes the order of $a\in G$ relative to $H$
(that is, the sma\-llest positive integer $m$ such that $a^m\in
H$). Study the connections between $\psi(G)$ and the collection
$(\psi_H(G))_{H\leq G}$, as well as the minimum/maximun of
$\{\psi_H(G)\hspace{1mm}\mid\hspace{1mm} H\leq G, \mid
H\mid\,=n\}$, where $n\in \mathbb{N}^{*}$ is fixed.
\bigskip

\bigskip\noindent{\bf Acknowledgements.} The authors are grateful to the reviewers for
their remarks which improve the previous version of the paper.

\vspace*{5ex} \small

\begin{minipage}[t]{5cm}
Marius T\u arn\u auceanu \\
Faculty of  Mathematics \\
``Al.I. Cuza'' University \\
Ia\c si, Romania \\
e-mail: {\tt tarnauc@uaic.ro}
\end{minipage}
\hfill
\begin{minipage}[t]{7cm}
Dan Gregorian Fodor \\
Faculty of Mathematics \\
``Al.I. Cuza'' University \\
Ia\c si, Romania \\
e-mail: {\tt dan.ms.chaos@gmail.com}
\end{minipage}

\end{document}